\documentclass{noticesOne}
\usepackage{amsfonts,amssymb,amsmath}
\usepackage{graphicx}
\usepackage{booktabs}

\newcommand{\R}{\mathbb{R}}

\title{From Coefficients to Distributions:\\
De~Moivre and the Operational View of Probability}

\author{
  Rodrigo Labouriau
  \affil{Rodrigo Labouriau, Depatment of Mathematics, 
  Aarhus University. E-mail: \texttt{rodrigo.laboriau@rlstatlab.com}.}
}
\date{Spring 2026}

\begin{document}

\maketitle

\vspace{-1cm}

\section*{}

How do we \emph{know} a probability law?
The question seems settled: a probability measure is a
countably additive set function on a $\sigma$-algebra,
and we know it when we can evaluate it on measurable sets.
Yet for more than two centuries before Kolmogorov's
axiomatisation, probabilists extracted detailed quantitative
information from probability laws without any measure-theoretic
machinery.
They did so by probing distributions with algebraic
objects---coefficients of polynomials, values of generating
functions, moments---and reading off probabilistic information
from the answers.

This article traces a single structural idea through
three hundred years of mathematics:
\emph{a probability law is known through what it does to
a family of test objects.}
The idea begins with De~Moivre's
coefficient extraction in the 1730s,
passes through the generating functions of Euler and Laplace,
the characteristic functions of L\'evy and Cram\'er,
and arrives at the distributional pairings
$\langle T, \varphi\rangle$ of modern functional analysis.
At each stage, the ``probes'' become more
general, and the class of laws that can be studied grows wider.
The mathematical distance between De~Moivre's binomial
coefficients and tempered distributions is enormous.
The structural kinship, however, is unmistakable.

Along the way, we encounter a historical injustice worth
correcting.
The normal distribution is universally called ``Gaussian,''
but the first systematic derivation of the normal curve---its shape,
its normalising constant, and its tail probabilities computed
to six decimal places---belongs to Abraham de~Moivre,
who obtained it in 1733, more than seventy years before Gauss.

The note grew out of a conversation with
J\o rgen Hoffmann-J\o rgensen about the pre-measure-theoretic
roots of probability.
Reading De~Moivre's original text afterwards, I was struck
by how naturally his methodology maps onto ideas that have
recently resurfaced in a distributional approach to
statistics \cite{LabouriauA1, LabouriauTransversality}.
\begin{figure}[ht]
\centering
\includegraphics[width=0.95\textwidth]{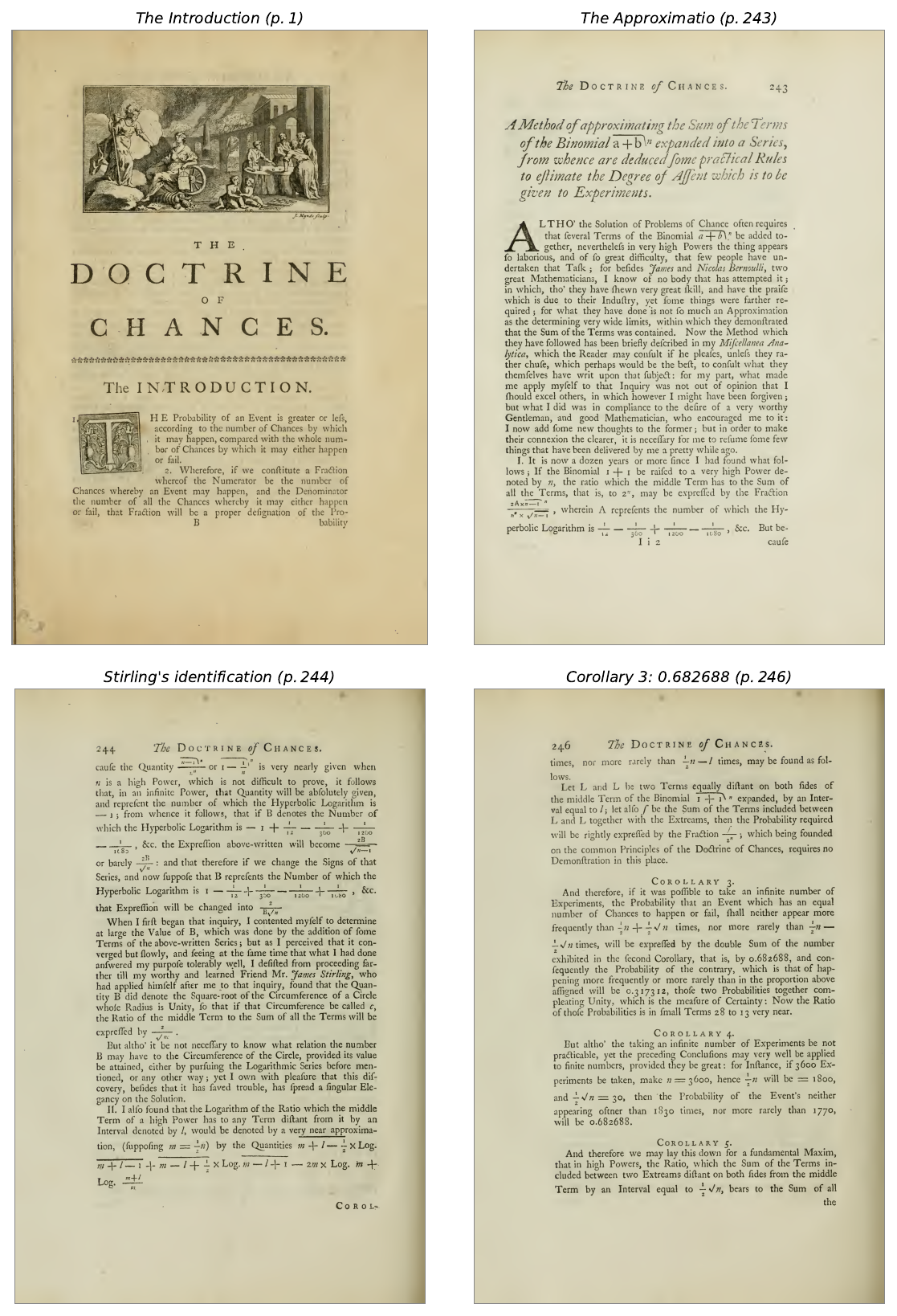}
\caption{Four pages from De~Moivre's \emph{Doctrine of Chances}
(1738).
\emph{Top left:} the opening of the book, with the definition
of probability as a fraction.
\emph{Top right:} the title of the \emph{Approximatio},
the pamphlet containing the first derivation of the normal curve.
\emph{Bottom left:} the passage where Stirling identifies
$B = \sqrt{2\pi}$, completing the normalising constant.
\emph{Bottom right:} Corollary~3, computing the $1$-$\sigma$
probability as $0.682688$.}
\label{fig:facsimile}
\end{figure}
\clearpage 

\section*{The pamphlet}

On November~12, 1733, Abraham de~Moivre (1667--1754)
circulated privately among friends a seven-page Latin pamphlet
entitled \emph{Approximatio ad Summam Terminorum Binomii
$\overline{a+b}\,^n$ in Seriem expansi}.
Only two copies are known to survive.
De~Moivre later translated and expanded the text into English,
incorporating it into the second edition of his
\emph{Doctrine of Chances} \cite{DeMoivre1738}
(pages 235--243; see also Walker's annotated
reprint \cite{Walker1934}).
The pamphlet gave the first derivation of the normal
curve, the first method for computing normal tail probabilities,
and the first recognition that the standard deviation is the
natural unit for measuring deviations from the mean.

De~Moivre begins with the symmetric case $p=q=\tfrac{1}{2}$.
The binomial $(1+1)^n$ has $2^n$ as its total sum, and the
middle coefficient $\binom{n}{n/2}$ is the largest term.
His goal is to approximate the ratio of any term to this
middle term when $n$ is large.
Using what we now call Stirling's series for $\log(n!)$,
he shows that
$$
\frac{\binom{n}{n/2+l}}{\binom{n}{n/2}}
\;\approx\;
\exp\!\left(-\frac{2l^2}{n}\right)
$$
when $l$ is small compared to~$n$.
This is the Gaussian exponent.
After the substitution $x = 2l/\sqrt{n}$, it becomes the
familiar $e^{-x^2/2}$.

De~Moivre had the shape of the curve, but not the
normalising constant.
His asymptotic expansion involved an unknown quantity~$B$:
the ratio of the middle term to the sum of all terms was
$2B/\sqrt{nc}$, where $c$ denotes the circumference of the
unit circle.
At this point occurs one of the most appealing episodes
of mathematical collaboration in the eighteenth century.
De~Moivre writes that he
``desisted from proceeding farther, till my worthy
and learned Friend Mr.\ \emph{James Stirling},
who had applied himself after me to that inquiry,
found that the Quantity~$B$ did denote the
Square-root of the Circumference of a Circle whose
Radius is Unity'' \cite{Walker1934}.
Stirling had identified $B = \sqrt{2\pi}$.
With this, the approximation becomes complete:
$$
\frac{\binom{n}{n/2+l}}{2^n}
\;\approx\;
\frac{1}{\sqrt{\pi n/2}}\;
\exp\!\left(-\frac{2l^2}{n}\right),
$$
which is the density of $\mathcal{N}(0,\,n/4)$ evaluated
at the point $n/2 + l$.

But De~Moivre did not stop at the formula.
He computed tail probabilities by numerical quadrature,
using Newton--Cotes rules (which he calls
``the Artifice of Mechanic Quadratures,
first invented by Sir~\emph{Isaac Newton}'').
His results are remarkable.

\medskip
\begin{center}
\begin{tabular}{ccc}
\toprule
Interval & De~Moivre & Modern \\
\midrule
$\tfrac{n}{2} \pm \tfrac{1}{2}\sqrt{n}$ \enspace ($\pm 1\sigma$)
  & $0.682688$ & $0.682689$ \\[3pt]
$\tfrac{n}{2} \pm \sqrt{n}$ \enspace ($\pm 2\sigma$)
  & $0.95428$ & $0.95450$ \\[3pt]
$\tfrac{n}{2} \pm \tfrac{3}{2}\sqrt{n}$ \enspace ($\pm 3\sigma$)
  & $0.99874$ & $0.99730$ \\
\bottomrule
\end{tabular}
\end{center}

\smallskip\noindent
The first line is the probability that
the count of successes in $n$ fair trials falls within
one standard deviation of the mean---what any statistics
student today would look up in a normal table---computed
to six correct decimal places in~1733.
Gauss would not introduce the normal distribution
in the context of error theory until 1809 \cite{Gauss1809};
Laplace's general version of the central limit theorem
appeared in 1812~\cite{Laplace1812}.

\begin{figure}[ht]
\centering
\includegraphics[width=\textwidth]{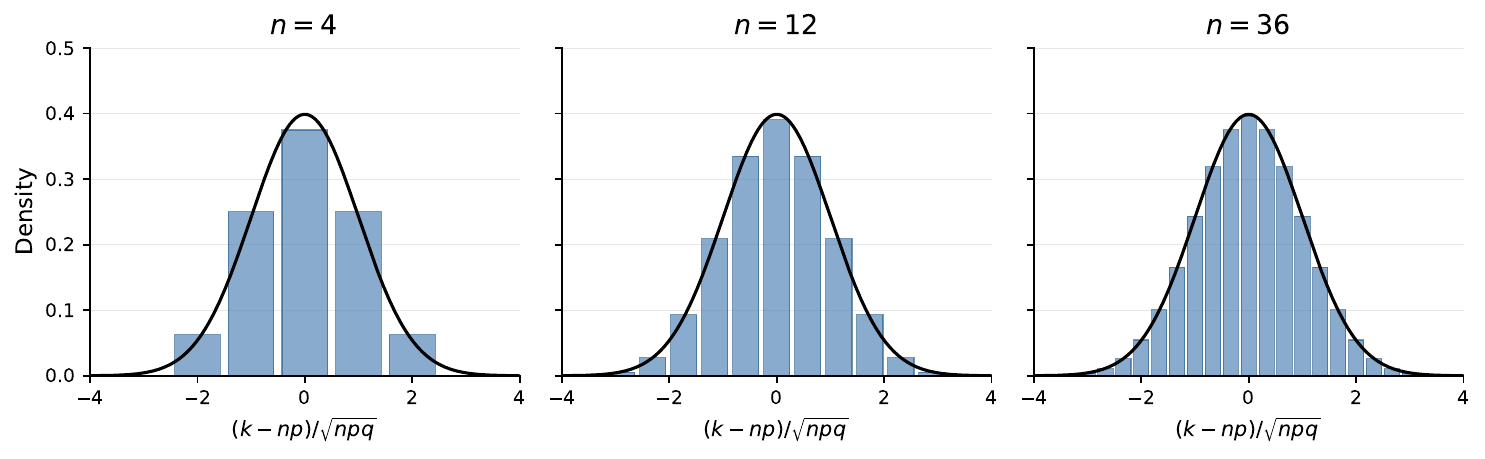}
\caption{De~Moivre's approximation visualised.
The standardised binomial probability function
(bars) converges to the standard normal density
(curve) as $n$ increases ($p = \tfrac{1}{2}$).}
\label{fig:binomial}
\end{figure}

De~Moivre then extends his analysis to the
general binomial $(a+b)^n$ with $a \neq b$.
In his Corollary~9,
the log-ratio becomes $-(a+b)^2 l^2/(2abn)$---again
a Gaussian exponent, now with variance $abn/(a+b)^2$,
which is $npq$ in modern notation.
The universality of the Gaussian shape across the
parameter space of Bernoulli models was already visible
to De~Moivre.

As Stigler \cite{Stigler1986} has documented in detail,
De~Moivre's priority is beyond question.
The normal distribution might more justly be called the
\emph{De~Moivre distribution}, or perhaps the
\emph{De~Moivre--Stirling distribution}.

\section*{Coefficients as probes}

What is De~Moivre actually doing in the \emph{Approximatio}?
He is not constructing a sample space.
He is not defining a $\sigma$-algebra or a measure.
He is extracting information about a probability law by
evaluating it against specific \emph{probes}:
the coefficient of the $k$-th term in the expansion of
$(p+q)^n$ gives $\Pr(X=k)$;
the sum of coefficients over an interval gives the probability
of that interval.

In modern language, the expansion
$$
(p+q)^n
\;=\;
\sum_{k=0}^{n} \binom{n}{k}\,p^k\,q^{n-k}
$$
defines the probability generating function evaluated at $z=1$.
At $z=e^{it}$, it becomes the characteristic function
$(pe^{it}+q)^n$.
And the characteristic function is itself a pairing:
if $\mu_n$ denotes the binomial measure, then
$(pe^{it}+q)^n = \int e^{itx}\,d\mu_n(x)$.

What De~Moivre computes in his Corollaries~3--5 is
even more directly a pairing.
He evaluates
$$
\sum_{k=0}^{n} \binom{n}{k}\,p^k\,q^{n-k}\;
\mathbf{1}_{[a,b]}\!\left(\frac{k-np}{\sqrt{npq}}\right)
$$
---the action of the standardised binomial distribution
on the indicator of an interval---and shows that
it converges to the Gaussian integral
$$
(2\pi)^{-1/2}\int_a^b e^{-x^2/2}\,dx.
$$

The structure is recognisably that of weak convergence,
even though the modern language of test-function topologies
lay far in the future.

\section*{De~Moivre's theorem in distributional language}

The observation above can be made precise.
Let $X_1, X_2, \ldots$ be independent Bernoulli random variables
with success probability $p$ and let $q = 1-p$.
The binomial sum $S_n = X_1 + \cdots + X_n$, after centring and
scaling, defines a standardised distribution $T_n$ that acts on
test functions by
$$
\langle T_n, \varphi\rangle
\;=\;
\sum_{k=0}^{n} \binom{n}{k}\,p^k\,q^{n-k}\;
\varphi\!\left(\frac{k-np}{\sqrt{npq}}\right).
$$

\noindent
\textbf{Proposition}
(De~Moivre--Laplace, distributional form)\textbf{.}
\textit{For every Schwartz function
$\varphi \in \mathcal{S}(\R)$,}
$$
\langle T_n, \varphi\rangle
\;\longrightarrow\;
\frac{1}{\sqrt{2\pi}}
\int_{\R} \varphi(x)\,e^{-x^2/2}\,dx
\;=:\;
\langle T_{\mathcal{N}},\,\varphi\rangle
\qquad
\textit{as } n \to \infty.
$$

\smallskip\noindent
\textit{Sketch.}
Stirling's formula gives the local approximation
$\binom{n}{n/2+l}\big/\binom{n}{n/2}\!\approx\!e^{-2l^2/n}$
for the individual terms.
Since $\varphi$ is rapidly decreasing, the tails of the sum are
negligible, and the local approximation controls the bulk.
The convergence is uniform over bounded subsets of
$\mathcal{S}(\R)$, so $T_n \to T_{\mathcal{N}}$ in
$\mathcal{S}'(\R)$.

\medskip
The proposition is a standard result, but its formulation
is worth pausing over.
De~Moivre himself proved the special case
$\varphi = \mathbf{1}_{[a,b]}$:
his Corollaries~3--5 compute the pairing
$\langle T_n, \mathbf{1}_{[-s,s]}\rangle$ for
$s = 1, 2, 3$ (in units of $\sigma$) and verify that
the result converges to the Gaussian integral.
The distributional formulation replaces indicator functions
by arbitrary Schwartz functions---a vast enlargement of the
probe family---but the conceptual structure is unchanged:
convergence is expressed as convergence of pairings against
probes.

\section*{Enlarging the probe family}

The history of probability after De~Moivre can be read,
in part, as a progressive enlargement of the family of probes
used to interrogate a probability law.

De~Moivre probes with \emph{coefficients}:
$\Pr(X=k)$ is extracted from the $k$-th term of a polynomial
expansion.
Euler and Laplace develop \emph{generating functions}:
the entire sequence of probabilities is encoded in
$(pz+q)^n$, and information is extracted by algebraic
operations on~$z$.
Chebyshev, Markov, and Stieltjes study the \emph{moment
sequence} $m_r = \sum k^r\,\Pr(X=k)$:
each moment is the pairing of the distribution against
the monomial probe $x^r$.
L\'evy and Cram\'er introduce the \emph{characteristic
function} $\phi(t) = \int e^{itx}\,d\mu(x)$:
the probe is now the exponential $e^{itx}$,
and the full distribution is recovered from the
function $t \mapsto \phi(t)$ by Fourier inversion.

At each step, the probes become more flexible and the
class of distributions that can be studied grows.
Coefficients work for discrete distributions on finitely
many points.
Generating functions handle sequences.
Moments characterise many distributions, but not all---the
classical moment problem has indeterminate solutions.
Characteristic functions determine \emph{every} probability
measure, but individual moments may fail to exist
(think of the Cauchy distribution, where $\mathbb{E}[|X|]$
is already infinite).

There is a natural further step.
Instead of probing a distribution with a single monomial
$x^r$ or a single exponential $e^{itx}$,
one can probe it with an arbitrary
\emph{test function} $\varphi$ from a suitable class.
If $T$ is a tempered distribution
(in the sense of Schwartz~\cite{Schwartz1966})
and $\varphi$ belongs to the Schwartz space
$\mathcal{S}(\R)$---the space of smooth, rapidly decreasing
functions---then the pairing
$\langle T, \varphi\rangle$ is always well defined.
The product $x^r\varphi(x)$ again belongs to
$\mathcal{S}(\R)$ (Schwartz space is closed under
multiplication by polynomials), so
$\langle T, x^r\varphi\rangle$ exists for every~$r$:
all ``moments'' are finite, even for distributions like the
Cauchy where the classical moments diverge.
The pairing $\langle T, e^{itx}\varphi(x)\rangle$
likewise exists and defines a ``regularised characteristic
function.''

The chain of generalisations thus runs:
\begin{eqnarray}\nonumber
\underbrace{\textstyle\binom{n}{k}p^kq^{n-k}}_{\text{coefficient}}
 & \;\longrightarrow\; &
\underbrace{\textstyle\sum k^r \Pr(X\!=\!k)}_{\text{moment}}
\;\longrightarrow\;
\underbrace{\textstyle\int e^{itx}\,d\mu}_{\text{characteristic\ function}}
\;\longrightarrow\;
\underbrace{\langle T,\,\varphi\rangle}_{\text{pairing}}.
\end{eqnarray}
At each stage, the probe is more general, the class of
accessible distributions is larger, and the extraction of
probabilistic information is more robust.
In recent work \cite{LabouriauA1},
these pairings are used systematically to define
``weak moments,'' ``weak cumulants,'' and
a ``weak characteristic function'' for distribution--kernel
pairs $(T,\varphi)$, recovering and extending
the classical moment and cumulant calculus in a setting
where the underlying distribution need not possess
a density or finite moments of any order.

\section*{Two ways of knowing a probability law}

The Kolmogorov axiomatisation of probability~\cite{Kolmogorov1933}
proceeds along a definite chain:
\begin{eqnarray}\nonumber
\text{sample space}
& \;\rightarrow\; &
\sigma\text{-algebra} 
 \;\rightarrow\
\text{probability measure}
\;\rightarrow\;
\text{expectation}.
\end{eqnarray}
This foundational sequence is so deeply embedded in modern
teaching that it is easy to forget it is a choice, not a
logical necessity.
There is another operational viewpoint:
\begin{eqnarray}\nonumber
\text{test functions}
& \;\rightarrow\; &
\text{distribution--kernel pairs}
\;\rightarrow\;
\text{weak expectations}.
\end{eqnarray}
This viewpoint does not reject the first---every
probability measure defines a tempered distribution,
and the distributional expectation
$\langle T, \varphi\rangle$ agrees with the measure-theoretic
one whenever both are defined.
But it is available in situations where the
measure-theoretic route is not:
when the ``density'' is too singular to be a
Radon--Nikodym derivative,
when moments of a given order do not exist in the classical
sense, or when the distribution is specified operationally
(through a transform or a characteristic function) rather
than constructively (through a density with respect to
a base measure).

It is this operational viewpoint that is structurally
closer to the pre-Kolmogorov tradition.
De~Moivre, Euler, Laplace, and Chebyshev all
extracted information from distributions by algebraic
manipulation, asymptotic expansion, and coefficient comparison, all
without measure-theoretic foundations.
The distributional approach returns to that operational
philosophy, now equipped with the full power of modern
functional analysis.\par

\section*{A puzzle and its geometric resolution}

There is a classical puzzle that connects to this story.
The moment problem asks:
given a sequence $(m_0, m_1, m_2, \ldots)$,
is there a unique probability measure with those moments?
For most distributions encountered in practice,
the answer is yes---the distribution is
\emph{M-determinate}.
Yet the theory of Stieltjes classes shows that
M-indeterminate distributions exist in abundance.
Why, then, are they so rarely encountered
in parametric statistical models?

The answer turns out to be geometric.
In recent work (\cite{LabouriauTransversality}, Example~8.2),
it is shown that the set of M-indeterminate distributions
has infinite codimension in the space of all distributions,
so that a ``generic'' parametric family avoids it entirely.
This is a transversality phenomenon:
the parametric family is a finite-dimensional submanifold,
and the M-indeterminate set is a thin stratum of infinite
codimension; they miss each other for dimensional reasons.
The same geometric mechanism explains other statistical
abnormalities: non-identifiability corresponds to
self-intersections of the feature map, singular Fisher
information to rank deficiency, and higher-order
instabilities to degeneracies of higher jets.
In each case the offending set has high codimension,
and transversality ensures that generic parametric
families avoid it.

De~Moivre's binomial family is trivially M-determinate,
being supported on finitely many points.
But its Gaussian limit is M-determinate for a deeper
reason: exponential tail decay ensures that the moment
generating function converges in a neighbourhood of the
origin, which is Carleman's sufficient condition
(\cite{LabouriauTransversality}, Appendix~A).
The distributional machinery captures this naturally
by embedding the Gaussian into a tilted measure whose
moments grow slowly enough to guarantee determinacy.
The passage from the combinatorial finiteness of
De~Moivre's binomial to the analytic structure of the
Gaussian limit illustrates, in miniature, the passage
from pre-modern operational probability to the modern
distributional framework.

\section*{Coda}

De~Moivre did not possess measure theory, tempered
distributions, Schwartz spaces, or weak convergence.
But the structural architecture of his work is visible
in retrospect:
probability laws characterised through algebraic expansions;
asymptotic approximation by a universal limiting shape;
information extracted through the evaluation of distributions
against probes.

De~Moivre probed the binomial law with
coefficients and indicator functions of intervals.
The distributional approach probes tempered distributions
with Schwartz functions.
The mathematical distance is enormous;
the structural kinship is unmistakable.

 \smallskip
\noindent\textsc{Acknowledgement.}
I am grateful to J\o rgen Hoffmann-J\o rgensen
for the conversation that prompted a re-reading of
De~Moivre's text.


 \newpage
\appendix
\section*{Appendix A: Proof of the Proposition}
\label{app:proof}

We prove that $\langle T_n, \varphi\rangle \to
\langle T_{\mathcal{N}}, \varphi\rangle$ for every
$\varphi \in \mathcal{S}(\R)$, where
$$
\langle T_n, \varphi\rangle
\;=\;
\sum_{k=0}^{n} \binom{n}{k}\,p^k\,q^{n-k}\;
\varphi\!\left(\frac{k-np}{\sqrt{npq}}\right).
$$
Write $x_k = (k - np)/\sqrt{npq}$ for the standardised
argument and $w_k = \binom{n}{k}\,p^k\,q^{n-k}$ for the
binomial weight.

\medskip\noindent
\emph{Local Gaussian approximation:}
For $|x_k|$ bounded (say $|x_k| \le M$), Stirling's formula
yields
$$
w_k
\;=\;
\frac{1}{\sqrt{2\pi npq}}\,
\exp\!\left(-\frac{x_k^2}{2}\right)
\bigl(1 + O(n^{-1/2})\bigr),
$$
where the $O$-constant depends only on $M$ and $p$.
This is De~Moivre's own computation, extended from
$p=\tfrac{1}{2}$ to general $p$ in his Lemma~3.

\medskip\noindent
\emph{Tail estimate:}
Since $\varphi \in \mathcal{S}(\R)$, for every $N > 0$
there exists $C_N$ such that $|\varphi(x)| \le C_N(1+|x|)^{-N}$.
The Chebyshev bound gives
$$
\sum_{|x_k| > M} w_k \;=\;
P\bigl(|S_n - np| > M\sqrt{npq}\bigr)
\;\le\; \frac{1}{M^2}.
$$
The contribution of the tail terms to
$\langle T_n, \varphi\rangle$ is therefore bounded by
$C_N\,M^{-2}$, uniformly in~$n$.

\medskip\noindent
\emph{Riemann sum approximation:}
On the bulk region $|x_k| \le M$, the spacing between
consecutive $x_k$ values is $\Delta x = 1/\sqrt{npq}$.
Using the local Gaussian approximation,
\begin{align*}
\sum_{|x_k|\le M} w_k\,\varphi(x_k)
&\;=\;
\sum_{|x_k|\le M}
\frac{1}{\sqrt{2\pi}}\,e^{-x_k^2/2}\,
\varphi(x_k)\,\Delta x
\;+\; O(n^{-1/2})\\
&\;\longrightarrow\;
\frac{1}{\sqrt{2\pi}}\int_{-M}^{M}
\varphi(x)\,e^{-x^2/2}\,dx.
\end{align*}

\medskip\noindent
\emph{Conclusion:}
Letting $M \to \infty$ (after $n \to \infty$),
the Gaussian tail integral vanishes because
$\varphi\,e^{-x^2/2} \in L^1(\R)$.
Combining the tail estimate and the
Riemann sum approximation:
$$
\langle T_n, \varphi\rangle
\;\longrightarrow\;
\frac{1}{\sqrt{2\pi}}\int_{\R}
\varphi(x)\,e^{-x^2/2}\,dx
\;=\;
\langle T_{\mathcal{N}}, \varphi\rangle.
$$
The convergence is uniform over bounded subsets of
$\mathcal{S}(\R)$ because the seminorm bounds on
$\varphi$ enter only through the constants $C_N$
and the smoothness of the Riemann-sum
approximation.  Hence $T_n \to T_{\mathcal{N}}$
in $\mathcal{S}'(\R)$.
\hfill$\square$

\newpage
\section*{Appendix B: Guide to the \emph{Approximatio}}
\label{app:guide}

The following table locates the key results of De~Moivre's
\emph{Approximatio} within the Walker
reprint~\cite{Walker1934}.

\smallskip
{\small
\begin{tabular}{p{3.3cm}p{3.7cm}p{0.6cm}}
\toprule
\textbf{Text} & \textbf{Content} & \textbf{p.} \\
\midrule
Opening
  & Stirling series for $\log(n!)$;
    ratio $2B/\sqrt{nc}$
  & 75 \\[2pt]
Stirling's identification
  & $B = \sqrt{2\pi}$
  & 76 \\[2pt]
Corollary 1
  & Gaussian exponent:
    $\log\text{-ratio} = -2l^2/n$
  & 77 \\[2pt]
Corollary 2
  & Series for the
    binomial sum
  & 77 \\[2pt]
Lemma
  & CLT statement
    for $p=\tfrac{1}{2}$
  & 78 \\[2pt]
Corollary 3
  & $1$-$\sigma$: $0.682688$
  & 78 \\[2pt]
Corollary 5
  & $2$- and $3$-$\sigma$
    probabilities
  & 79 \\[2pt]
Corollary 6
  & Newton--Cotes
    quadrature
  & 80 \\[2pt]
Lemma 3
  & Extension to
    $p \neq \tfrac{1}{2}$
  & 82 \\[2pt]
Corollary 8
  & Normalising const.\
    $(a\!+\!b)/\sqrt{abnc}$
  & 82 \\[2pt]
Corollary 9
  & General Gaussian
    exponent
  & 83 \\[2pt]
Corollary 10
  & Regularity from
    chance; Design
  & 83 \\
\bottomrule
\end{tabular}
}


\begin{thebibliography}{9}

\bibitem{DeMoivre1738}
A.~de~Moivre,
\emph{The Doctrine of Chances},
2nd ed., Woodfall, London, 1738.

\bibitem{Gauss1809}
C.~F.~Gauss,
\emph{Theoria motus corporum coelestium},
Perthes et Besser, Hamburg, 1809.

\bibitem{Kolmogorov1933}
A.~N.~Kolmogorov,
\emph{Grundbegriffe der Wahrscheinlichkeitsrechnung},
Springer, Berlin, 1933.
MR0494348

\bibitem{LabouriauA1}
R.~Labouriau,
Distributional statistical models: weak moments, cumulants,
and a central limit theorem,
\emph{arXiv}:2604.20634 [math.PR], 2026.

\bibitem{LabouriauTransversality}
R.~Labouriau,
Transversality and geometric regularisation in distributional
statistical models.
\emph{arXiv}:2605.04536 [math.ST] 2026.

\bibitem{Laplace1812}
P.-S.~Laplace,
\emph{Th\'eorie analytique des probabilit\'es},
Courcier, Paris, 1812.

\bibitem{Schwartz1966}
L.~Schwartz,
\emph{Th\'eorie des distributions},
2nd ed., Hermann, Paris, 1966.
MR0209834

\bibitem{Stigler1986}
S.~M.~Stigler,
\emph{The History of Statistics: The Measurement of
Uncertainty before 1900},
Harvard Univ.\ Press, Cambridge, MA, 1986.
MR0852410

\bibitem{Walker1934}
H.~M.~Walker,
De~Moivre on the law of normal probability,
in D.~E.~Smith (ed.),
\emph{A Source Book in Mathematics},
McGraw-Hill, New York, 1929, pp.~566--575.

\end{thebibliography}
\end{document}